\documentclass{amsart}
\usepackage[utf8]{inputenc}
\usepackage{amsmath}
\usepackage{amsfonts}
\usepackage{amssymb}
\usepackage{enumitem,pstricks,xy,pst-node}
\xyoption{all}
\usepackage{amsthm}

\begin{document}


\def\cV{{\mathcal V}}\def\cM{{\mathcal M}}
\def\PP{\mathbb{P}}\def\AA{\mathbb{A}}\def\RR{\mathbb{R}}

\title[On a class of spaces of skew-symmetric forms]{On a class of spaces of skew-symmetric forms related to Hamiltonian systems of conservation laws}

\author{E.V. Ferapontov,  L. Manivel }


\date{}

\begin{abstract}
It was shown in \cite{FPV} that the  classification of $n$-component systems of conservation laws possessing a third-order Hamiltonian structure reduces to the following algebraic problem:
  classify $n$-planes $H$ in $\wedge^2(V_{n+2})$ such that  the induced map 
$Sym^2H\longrightarrow  \wedge^4V_{n+2}$ has 1-dimensional kernel  generated by a non-degenerate quadratic form on $H^*$.

This problem is trivial for $n=2, 3$ and apparently wild for $n\geq 5$. In this paper we address the most interesting borderline case $n=4$. 

We prove that the variety $\cV$ parametrizing those 4-planes $H$ is an irreducible 38-dimensional  $PGL(V_6)$-invariant subvariety
of the Grassmannian $G(4,  \wedge^2V_6)$. With every $H\in\cV$ we associate  a 
{\it characteristic} cubic surface $S_H\subset \PP H$, the locus of rank 4 two-forms  in $ H$.
We  demonstrate that the induced characteristic map $\sigma: \cV / PGL(V_6) \dashrightarrow \cM_c,$  where $\cM_c$ denotes the moduli space of cubic surfaces in $\PP ^3$, is dominant, hence generically finite. 
A complete classification of 4-planes $H\in\cV$ with the reducible characteristic  surface $S_H$ is given.

\end{abstract}



\maketitle

\setcounter{page}{1}

\newtheorem{theorem}{Theorem}[section]
\newtheorem{lemma}[theorem]{Lemma}
\newtheorem{proposition}[theorem]{Proposition}
\newtheorem{corollary}[theorem]{Corollary}
\newtheorem{definition}[theorem]{Definition}

\def\PP{\mathbb{P}}\def\AA{\mathbb{A}}\def\RR{\mathbb{R}}
\def\CC{\mathbb{C}}\def\HH{\mathbb{H}}\def\OO{\mathbb{O}}
\def\ZZ{\mathbb{Z}}\def\QQ{\mathbb{Q}}
\def\cA{{\mathcal A}}
\def\cC{{\mathcal C}}
\def\cV{{\mathcal V}}\def\cM{{\mathcal M}}

\def\cW{{\mathcal W}}\def\cT{{\mathcal T}}
\def\ra{\rightarrow}\def\lra{\longrightarrow}
\def\fg{{\mathfrak g}}\def\fe{{\mathfrak e}}\def\fsl{{\mathfrak{sl}}}\def\fso{{\mathfrak{so}}}
\def\s{{\sigma}}\def\t{{\tau}}\def\iot{{\iota}}
\def\GW{{Gromov-Witten}}

\tableofcontents

\newpage

\section{Introduction}

Our  problem originates from the geometric theory of systems of conservation laws,
\begin{equation} \label{V} 
u^i_t = [V^i(\mathbf{u})]_x,
\end{equation}
$i=1, \dots, n$. PDEs of this type  appear in a wide range of applications in continuum mechanics and mathematical
physics, see e.g. \cite{Ser}.
A geometric counterpart of system (\ref{V})
is a line congruence ($n$-parameter family of lines in projective space $\mathbb{P}^{n+1}$) specified by the equations 
\begin{equation}
  y^i=u^i y^{n+1}+V^i({\bf u})y^{n+2},
  \label{cong}
\end{equation}
here $y^i$ are the  homogeneous coordinates in $\mathbb{P}^{n+1}$.
In the case $n=2$ we obtain a two-parameter
 family, or a classical congruence of lines in $\mathbb{P}^3$. Since 19th century the theory of
 congruences has been one of the most popular chapters of projective  differential
 geometry.
It was observed in \cite{AF1, AF2} that all standard concepts of the
theory of conservation laws such as rarefaction curves, shock curves, linear
degeneracy, reciprocal transformations, etc, acquire a simple intuitive interpretation in
the language of the projective theory of congruences.   Algebro-geometric aspects  of the correspondence (\ref{V}) $\leftrightarrow$ (\ref{cong}) were thoroughly investigated in \cite{Mez1, Mez2}.

Particularly interesting examples of systems (\ref{V}) arise in the context of equations of  associativity of 2D topological field theory (WDVV equations) \cite{FN, pv}. Such systems can be represented in  Hamiltonian form
\begin{equation} \label{H}
u^i_t=P^{ij}\frac{\delta H} {\delta u^j}
\end{equation}
where $P^{ij}$ is a third-order Hamiltonian operator of special type 
and $H$ is a (necessarily nonlocal) Hamiltonian, see \cite{FPV} for further details. It was shown in \cite{FPV} that if system (\ref{V}) possesses  Hamiltonian representation (\ref{H}) then the associated congruence (\ref{cong}) must necessarily be linear, that is,  defined by $n$ linear equations in the  Pl\"ucker coordinates. Explicitly, congruence (\ref{cong}) must satisfy the relations
\begin{equation}
trYA^{i}=0,
\label{Y}
\end{equation}
$i=1, \dots, n$, where $Y$ is the $(n+2)\times (n+2)$ skew-symmetric matrix formed by $2\times 2$ minors (Pl\"ucker coordinates) of the $2\times (n+2)$ matrix
$$
\left(
  \begin{array}{ccccc}
    u^i & \dots &u^n&1&0\\
    V^i & \dots &V^n&0&1
  \end{array}
\right),
$$
and $A^i$ are constant $(n+2)\times (n+2)$ skew-symmetric matrices. Note that system (\ref{Y}) can be viewed as a linear system for the fluxes $V^i$ of system (\ref{V}). Furthermore, viewed as 2-forms,  $A^i$ must satisfy an additional relation of the form
\begin{equation}\label{varphi}
\varphi_{ij}A^i\wedge A^j=0
\end{equation}
where the matrix $\varphi$ is symmetric and non-degenerate. 

Introducing an $n$-plane $H=span \langle A^i \rangle$ in $\wedge^2(V_{n+2})$, interpreting relation (\ref{varphi}) as a kernel of the natural map
$Sym^2H\longrightarrow  \wedge^4V_{n+2}$, and $\varphi$ as a non-degenerate quadratic form on $H^*$, we arrive at the algebraic problem formulated above. In this paper we will concentrate on the particularly interesting case $n=4$.

\section{Problem and strategy}

\subsection{The problem}

As outlined in the introduction, 
we would like to classify the four-planes $H$ in $\wedge^2V_6$ such that
\begin{enumerate}
\item 
the induced map 
$Sym^2H\longrightarrow  \wedge^4V_6\simeq\wedge^2V_6^*$
has rank exactly nine, 
\item its kernel  is generated by a non-degenerate quadratic form on $H^*$. 
\end{enumerate}

The variety $\cV$ parametrizing those four-planes is a $PGL(V_6)$-invariant subvariety
of the Grassmannian $G(4,  \wedge^2V_6)$, whose dimension is $44$. The condition that 
$Sym^2H\longrightarrow  \wedge^4V_6$ be of rank at most nine describes a closed subvariety
$\cV^*$ whose expected codimension is $15-10+1=6$, so the expected dimension is $38=\dim
PGL(V_6)+3$. Then $\cV$ is defined by two open conditions in $\cV^*$, so it must be open in 
$\cV^*$ but not necessarily dense since it is not clear that $\cW$ is irreducible. More precisely,
$\cV$ is a union of dense open subsets of irreducible components of $\cV^*$. Moreover all these
components have dimension at least $38$. 

\subsection{Ranks} 

Elements of $\wedge^2V_6$ are skew-symmetric forms on $V_6^*$, for which we have the usual notion of rank. For convenience we will also 
use the following terminology.
\begin{definition}
Let $H$ be a subspace of $\wedge^2V_6$. Let $\theta$ be an element of
$Sym^2H$. 
\begin{enumerate} 
\item The {\it q-rank} of $\theta$ is the rank of the corresponding quadratic form on $H^*$;
\item The {\it rank} of $\theta$ is the rank of its image in $\wedge^4V_6\simeq \wedge^2V_6^*$, considered as a skew-symmetric form on $V_6$. 
\end{enumerate}
\end{definition}

For a given $\theta$, the q-rank depends on $H$, contrary to the 
rank. 

\medskip
Note that a four-plane $H$ in $\cV$ cannot contain any form $\omega$ of rank two, since such a form verifies $\omega\wedge\omega=0$, hence 
produces a degenerate element (of q-rank one) in the kernel of the
map $Sym^2H\longrightarrow  \wedge^4V_6$. As a consequence, $H$ 
contains elements of rank six: it was proved in \cite{mm} that the linear spaces whose non 
zero elements all have rank four have dimension at most three, and such three-planes 
were classified. 

\subsection{The birational involution}

Recall that the Pfaffian quadratic map
$$\begin{array}{rcl}
pf : \wedge^2V_6 & \longrightarrow & \wedge^4V_6 \\
  \omega & \mapsto & \omega\wedge\omega 
\end{array}$$
is birational, and essentially an involution. Indeed, if we identify $\wedge^4V_6$ to $\wedge^2V_6^*$, hence  $\wedge^4V_6^*$
to  $\wedge^2V_6$, we get a map
$$\begin{array}{rcl}
pf^* : \wedge^4V_6 & \longrightarrow & \wedge^2V_6 \\
  \omega^* & \mapsto & \omega^*\wedge\omega^*,
\end{array}$$
and $pf^*\circ pf(\omega)=Pf(\omega)\omega$, where $Pf(\omega)$ is essentially $\omega\wedge\omega \wedge\omega $, or the Pfaffian of
$\omega$ (which is defined only up to a non zero scalar). So if $\omega$ has rank six, we can recover it from $\omega\wedge\omega$ which has
also rank six. Note that if $\omega$ has rank four (resp. two), then $\omega\wedge\omega$ has rank two (resp. zero). In particular, $\omega\wedge\omega$ never has rank four.

\subsection{Irreducibility}

\begin{theorem}
$\cV$ is irreducible of dimension $38$.
\end{theorem}

\begin{proof}
Consider a four-plane $H$ that belongs to $\cV$, and denote by $q$ a generator of the kernel 
of the map $Sym^2H\ra\wedge^4V_6$. As a quadratic form on $H^*$, the tensor $q$ is 
non degenerate, hence identifies $H^*$ with $H$, yielding a quadratic form $q^*$ on $H$. 
As we observed, a generic $\omega$ in $H$ has rank six and q-rank four. Let $L\subset H$ 
be its orthogonal with respect to $q^*$. Then we can write, up to scalar, 
$$q=\omega^2-\theta, \qquad \mathrm{with} \quad \Omega \in Sym^2L.$$ 
In other words, $\Omega=\omega\wedge\omega$ belongs to $Sym^2L\simeq L\wedge L$
(the map $Sym^2L\ra L\wedge L$ is injective, hence an isomorphism, because its kernel 
is contained in the kernel of the map $Sym^2H\lra\wedge^4V_6$, which by hypothesis does
not contain any non zero element of q-rank three or less). Moreover, since $\omega$
has rank six, $\Omega$ also has rank six, and therefore $\Omega\wedge\Omega$ is a non
zero multiple of $\omega$. In particular we can recover $H$ from $L$ and $\Omega\in Sym^2L$. 

This suggests to consider the diagram 
\begin{equation*}
\label{kempfCollapsing}
\xymatrix @C=2pc @R=0.4pc{
\PP(Sym^2U) \ar[dd] & \rule{1pt}{0pt}\cW \ar@{_{(}->}[l] \ar[dr]^-{\pi} \ar[dd]^-{p_2} \\
& & \rule{3pt}{0pt}\hspace*{3mm}\cV\subset G(4,\wedge^2V_6),\\
G(3,\wedge^2V_6) & \rule{1pt}{0pt}\cT \ar@{_{(}->}[l] 
}
\end{equation*}
our notations being the following. First, $\cT$ is the open subset of $G(3,\wedge^2V_6)$ parametrizing the three-planes $L$ such that 
\begin{enumerate}
\item the map $Sym^2L\lra\wedge^4V_6$ is injective,
\item the general element in $L\wedge L$  has rank six.
\end{enumerate}
Second, $U$ denotes the tautological rank three vector bundle on $G(3,\wedge^2V_6)$ and 
$\cW$ is the open subset of $\PP(Sym^2U) $ parametrizing pairs $(L,[\Omega])$ such that 
\begin{enumerate}
\item $L$ belongs to $\cT$,
\item $\Omega\in Sym^2L\simeq L\wedge L$ has maximal rank and q-rank,
\item $(\omega\wedge L)\cap (L\wedge L)=0$, where $\omega:=\Omega\wedge\Omega$. 
\end{enumerate}
The map $\pi$ from $\cW$ to $\cV$ sends the pair $(L,[\Omega])$ to $H=\langle L, 
\omega\rangle$. This is well-defined because
\begin{itemize}
\item $\omega$ cannot belong to $L$: otherwise, we would have 
$\omega\wedge\omega=Pf(\Omega)\Omega$, yielding the relation $\omega^2=Pf(\Omega)\Omega$ 
in $Sym^2L$, and then $\Omega$ would have q-rank one, a contradiction; 
\item the kernel of the map $Sym^2H\ra \wedge^4V_6$ is generated by $q=\omega^2-Pf(\Omega)\Omega$,
which is non degenerate since $\Omega$ has q-rank three: indeed, if there was another element in 
this kernel, we could write it, after substracting a suitable multiple of $q$ if necessary, 
as $q'=\omega\lambda-\Phi$ for some $\lambda$ in $L$ and some $\Phi$ in $Sym^2L\simeq L\wedge L$; 
this would yield a relation $\omega\wedge\lambda=\Phi$ in $\wedge^4V_6$, contradicting the 
hypothesis that $(\omega\wedge L)\cap (L\wedge L)=0$. 
\end{itemize}

By the remarks at the beginning of this proof, $\pi$ is surjective, so the irreducibility
of $\cV$ will be a consequence of the irreducibility of $\cW$. But the latter is 
obviously an open subset of the projective bundle $\PP(Sym^2U)$, hence certainly 
irreducible, of dimension $3\times 12+5=41$. Finally, the fibers of $\pi$ are open
subsets in three-dimensional projective spaces, and therefore the dimension of $\cV$
is $41-3=38$. 
\end{proof}

\subsection{Characteristic surfaces} 

For $H\in\cV\subset G(4,\wedge^2V_6)$,  the locus of rank four two-forms defines a 
cubic surface $S_H\subset \PP H$, that we call the {\it characteristic surface}. 
This surface is the intersection of $\PP H$ with the Pfaffian hypersurface $Pf
\subset \PP \wedge^2V_6$. Recall that the singular locus of this hypersurface
is the Grassmannian $G(2,V_6)$, the locus of two-forms of rank two. Although 
$\PP H$ does not meet this Grassmannian when $H$ belongs to $\cV$, will we 
see later on that $S_H$ can be badly singular, and even reducible. 
Nevertheless, we have a natural map
$$\sigma : \cV / PGL(V_6) \dashrightarrow \cM_c,$$
where $\cM_c$ denotes the moduli space of cubic surfaces in $\PP ^3$. Both spaces
are four-dimensional, and we will see later on that $\sigma$ is dominant, hence generically
finite. 

Note that by construction,  $H\in\cV$ defines not only the cubic surface $S_H\subset \PP H$,
but also a Pfaffian representation of this surface. This is a classical topic. For example,
any cubic surface (possibly singular) admits a Pfaffian representation 
\cite[Proposition 7.6]{beauville}.  Moreover a generic such representation can be determined 
by five points in general position on the surface\cite[Example 7.4]{beauville}. A simple
algorithm for that is presented in \cite{tanturri}. 

The family of Pfaffian representations of a smooth cubic surface is five dimensional.
In comparison, there are only finitely many determinantal representations. They are  
in bijection with the $72$ linear systems of twisted cubics on the surface, or the 
$72$ sixes of non incident lines (this was already known in the 19th century, see 
\cite[Corollary 6.4]{beauville}).

\section{Semisimple three planes}

\subsection{A grading of $\fe_7$}
Although the classification of orbits in $G(4,  \wedge^2V_6)$ is not known (a priori a wild 
problem), the classification of orbits  in $G(3,  \wedge^2V_6)$ is known. This is because 
$V_3\otimes  \wedge^2V_6$ is part of a $\ZZ_3$-grading of $\fe_7$:
$$\fe_7\simeq \fsl(V_3)\times \fsl(V_6) \oplus V_3\otimes  \wedge^2V_6
 \oplus V_3^*\otimes  \wedge^2V_6^*.$$
Then it makes sense to talk about semisimple or nilpotent elements in $V_3\otimes  \wedge^2V_6$, 
by considering them as elements of $\fe_7$. The set of semisimple elements is the union of the 
Cartan subspaces, which are all equivalent under the action of $GL(V_3)\times GL(V_6)$. In order 
to exhibit one of those, we fix a basis $v_1, v_2, v_3$ of $V_3$, and a basis $e_1, \ldots , e_6$ 
of $V_6$. We let $e_{ij}=e_i\wedge e_j\in \wedge^2V_6$. 

\begin{proposition}
The subspace $C$ generated by the following three vectors is a Cartan subspace of $V_3\otimes  \wedge^2V_6$:

$$c_1 = v_1\otimes e_{12}+ v_2\otimes e_{34}+ v_3\otimes e_{56},$$
$$c_2 = v_1\otimes e_{36}+ v_2\otimes e_{52}+ v_3\otimes e_{14},$$
$$c_3 = v_1\otimes e_{54}+ v_2\otimes e_{16}+ v_3\otimes e_{32}.$$
\end{proposition}

We say that a three plane $L$ in $ \wedge^2V_6$ is semisimple if it admits a basis 
$\omega_1, \omega_2, \omega_3$ such that $ v_1\otimes \omega_1+ v_2\otimes  \omega_2
+ v_3\otimes \omega_3$ is semisimple in   $V_3\otimes  \wedge^2V_6$. By the previous
proposition, this is the case if and only if there exists a basis $e_1, \ldots , e_6$ 
of $V_6$, and coefficients $a,b,c$, not all zero, such that 

$$\omega_1 = ae_{12}+ be_{36}+ ce_{54},$$
$$\omega_2 = ae_{34}+ be_{52}+ ce_{16},$$
$$\omega_3 = ae_{56}+ be_{14}+ ce_{32}.$$

 
\subsection{Four planes containing a semisimple three plane}. We would like to understand
four planes $H\in \cV$ that contain a semisimple three plane $L$. Necessarily, the map $Sym^2L\longrightarrow 
 \wedge^4V_6$ must be injective. Note that 
$$\frac{1}{2}\omega_1^2 = bce_{3654}+ ace_{1254}+abe_{1236},\qquad -\omega_2\omega_3=a^2e_{3654}+ b^2e_{1254}+c^2e_{1236},$$
$$\frac{1}{2}\omega_2^2 = bce_{1652}+ ace_{1634}+abe_{3452},\qquad -\omega_1\omega_3=a^2e_{1652}+ b^2e_{1634}+c^2e_{3452},$$
$$\frac{1}{2}\omega_3^2 = bce_{1432}+ ace_{3256}+abe_{1456},\qquad -\omega_2\omega_3=a^2e_{1432}+ b^2e_{3256}+c^2e_{1456},$$

\noindent 
so this injectivity condition is equivalent to the fact that the vectors $(a^2, b^2, c^2)$ and $(bc, ac, ab)$ are not colinear. 

\smallskip Now we would like to complete $L$ with a vector $\omega_0\notin L$ such that the four plane $H$ they span belongs to $\cV$. We can choose $\omega_0$ such that it generates the
orthogonal to $L$ in $H$ with respect to its non degenerate quadratic form (indeed, 
generically we may suppose that the restriction of this quadratic form to $L$ is non 
degenerate). This implies that $\omega_0\wedge\omega_0$ belongs to $L\wedge L$.

Now observe that we can decompose $V_6=V_+\oplus V_-$, where 
$V_+=\langle e_1, e_3, e_5\rangle$ and $V_-=\langle e_2, e_4, e_6\rangle$. With respect to this decomposition, we have 
$$\wedge^2V_6=\wedge^2V_+\oplus (V_+\otimes V_-)\oplus \wedge^2V_-,$$ 
and $L\subset V_+\otimes V_-$. Let us decompose accordingly $\omega_0=\omega_++\Omega+\omega_-$. 
We get 
$$\frac{1}{2}\omega_0^2=\omega_+\Omega +\Big( \omega_+\omega_-+\frac{1}{2}\Omega^2 \Big) + \Omega\omega_-,$$
with respect to the decomposition 
$$\wedge^4V_6=(\wedge^3V_+\otimes V_-)\oplus (\wedge^2V_+\otimes \wedge^2V_-)
\oplus (V_-\otimes \wedge^3V_+).$$
Of course $L\wedge L$ is contained in $\wedge^2V_+\otimes \wedge^2V_-$ so we need in 
particular that 
$$\omega_+\Omega=\Omega\omega_-=0.$$
Generically, $\Omega$ has rank three, and then the previous equations imply 
$\omega_+=\omega_-=0$. This means that $H\subset V_+\otimes V_-$. Then $H$ belongs to 
$\cW$, since the image of $Sym^2H$ is 
contained in $\wedge^2V_+\otimes \wedge^2V_-$, which has dimension nine. One can 
check that a generic such $H$ does belong to  $\cV$. In fact it suffices to produce an 
explicit example; here is one:
$$\omega_1=e_{34}-e_{56}, $$
$$\omega_2=e_{52}-e_{14}, $$
$$\omega_3=e_{16}-e_{32}, $$
$$\omega_0=e_{12}+e_{14}+e_{32}+e_{36}+e_{54}-e_{56}. $$
The relation is 
$$\omega_0^2+\omega_1^2-\omega_2^2+\omega_3^2+\omega_1\omega_2+\omega_1\omega_3-\omega_2\omega_3=0,$$
which is non degenerate. 
 
 \medskip 
 When $H\subset V_+\otimes V_-$, we get a $3\times 3$ matrix of linear forms in four 
 variables, and the characteristic surface $S_H$ is given by the determinant of this matrix. Moreover the quadratic relation in $Sym^2H$ is orthogonal to the image of the map
 $$\wedge^2V_+^*\otimes \wedge^2V_-^*\longrightarrow Sym^2H^*$$
 given by the $2\times 2$ minors of our matrix. In the classical
 terminology, the quadratic relation in $Sym^2H$ is thus nothing
 else than the {\it Schur quadric} associated to the given determinantal 
 representation of the cubic surface $S_H$. When this surface is 
 smooth, it is known that the Schur quadric is unique and smooth (see \cite[9.1.3]{CAG}). 
 In particular this implies: 
 
 \begin{theorem} \label{moduli}
 The image of the characteristic map $\sigma : \cV/PGL(V_6)\lra \cM_c$ contains the open 
 subset parametrizing smooth cubic surfaces. In particular $\sigma$ is dominant, 
 hence generically finite. 
 \end{theorem}
 
 \noindent {\it Remark}. The stabilizer in $PGL(V_6)$ of a general point in $\cV$ is 
 the copy of $\CC^*$ given by the automorphisms of the form $tId_{V_+}+Id_{V_-}$. 
 This explains that the quotient $\cV/PGL(V_6)$ has dimension four, rather than $3=38-35$.
 
 \medskip
 We can also conclude from the previous discussion that in general, our special 
 Pfaffian representations are in fact determinantal. Since $\sigma$ is generically finite,
 this implies that $\cV/PGL(V_6)$ is birational to the moduli space of determinantal 
 representations of cubic surfaces. As a consequence, the degree of $\sigma$ must be $72$. 



\subsection{More examples}

Consider the three-plane $L$ generated by 
$$\begin{array}{rcl}
\omega_0 & = & e_0\wedge e_2 + e_1\wedge e_3, \\
\omega_1 & = & e_0\wedge e_4 + e_1\wedge e_5, \\
\omega_2 & = & e_0\wedge e_1+e_2\wedge e_5 + e_3\wedge e_4.
\end{array}$$
This is triple tritangent plane, the triple line being generic.
Let us complete it into a four-plane $H\in \cV$. For that we consider a general element 
of $L\wedge L$,
$$\Omega = xf_{23}+yf_{45}+z(f_{25}-f_{34})+a(f_{04}-
f_{15})+b(f_{02}-f_{13})+c(f_{01}+f_{25}+f_{34}).$$
Then we need to compute $\frac{1}{2}\Omega^2$. After simplyfying a little bit the result 
by adding a suitable element of $L$, we get 
$$\begin{array}{rcl}
\omega  & = & (xy-z^2)e_{01}+(yc+a^2)e_{23}+(xc+b^2)e_{45}
+ \\ & &+(yb+za)(e_{02}-e_{13})+(zb+xa)(e_{04}-e_{15})+
(zc+ab)(e_{34}-e_{25}).
\end{array}$$
Then we need to check that $H=\langle \omega, L\rangle$ belongs to $\cV$ when the parameters are general. We 
know that $\omega^2$ is a general element of $L\wedge L$, 
so the only thing we need to check is that the intersection of $L\wedge L$ with $\omega\wedge L$ is zero. Let $\Omega_1=f_{04}+f_{15}$, $\Omega_2=f_{02}+f_{13}$, $\Omega_3=f_{34}+f_{25}$. Note the space $\langle  
\Omega_1, \Omega_2, \Omega_3\rangle$ is transverse to 
$L\wedge L$. Moreover, modulo the latter we have
$$\begin{array}{rcl}
 \omega\omega_0 & =& -(zc+ab)\Omega_1-(xc+b^2)\Omega_2+(zb+xa)\Omega_3, \\
  \omega\omega_1 & =&
 -(yc+a^2)\Omega_1-(zc+ab)\Omega_2-(yb+za)\Omega_3, \\
  \omega\omega_2 & =&
 -(yb+za)\Omega_1-(zb+xa)\Omega_2-(zc+ab)\Omega_3.
 \end{array}$$
In general those three combinations of $\Omega_1, \Omega_2, \Omega_3$ are independent, and we are done. 

Note that the characteristic surface $S_H$ cannot be smooth. Indeed, by construction $L$ 
is a triple tritangent plane, and therefore $S_H$ must admit two singular points of 
type $A_2$, or even worse singularities (see e.g. \cite[9.2.2]{CAG} Dolgachev).

\medskip\noindent {\it Example 1}. Let $z,c$ be non zero and the other coefficients be zero, so (after factoring out a $z$) $\omega=-ze_{01}-c(e_{25}-e_{34})$. Here the (unique) quadratic relation is $$\frac{z}{2}\omega^2+c^2\omega\omega_2 =c(c^2-z^2)\omega_0\omega_1-\frac{zc^2}{2}\omega_2^2$$
and is non degenerate in general. The equation of the surface $S_H$ is 
$$Pf(t\omega+x_0\omega_0+x_1 \omega_1+x_2 \omega_2)=
zc^2t^3-c^2t^2x_2-ztx_2^2-2ctx_0x_1+x_2^3=0.$$
This surface has two singular points $[0,0,1,0]$ and $[0,1,0,0]$, each yielding an $A_2$ singularity.

\medskip\noindent {\it Example 2}. Let $b=y=1$ and the other coefficients be zero, so $\omega=e_{02}-e_{13}+e_{45}$. Here $\Omega=f_{02}-f_{13}+f_{45}$, so the (unique) quadratic relation is $$\frac{1}{2}\omega^2 =\omega_1\omega_2+\frac{1}{2}\omega_0^2$$
and is non degenerate. The equation of the characteristic surface $S_H$ is 
$$Pf(t\omega+x_0\omega_0+x_1 \omega_1+x_2 \omega_2)=
t^3-tx_0^2-2tx_1x_2+x_2^3=0.$$
This surface has a unique singularity at $[0,0,1,0]$ and this is an $A_4$ singularity.

\subsection{Special Pfaffian representations of  Cayley's ruled surface} 
It seems difficult to determine precisely what is the image of the characteristic map. 
We have just found some singular cubic surfaces with $A_2$ and even $A_4$ singularities,
and we will consider reducible surfaces in the next section. Among irreducible
surfaces, there are two types of non normal ones \cite[Theorem 9.2.1]{CAG}, and 
the most degenerate one is {\it Cayley's ruled surface}, of equation 
$$x_0^2x_2+x_1^2x_3=0.$$
In this section we discuss the special Pfaffian representations of  this surface. 

\medskip
So we consider $H=\langle \omega_0, \omega_1, \omega_2, \omega_3\rangle$ in $\cV$
and suppose that its characteristic surface $S_H$ has equation 
$$Pf(x_0\omega_0+x_1\omega_1+x_2\omega_2+x_3\omega_3)=x_0^2x_2+x_1^2x_3.$$
In particular, among the wedge products $\omega_i\omega_j\omega_k$, only 
$\omega_0^2\omega_2$ and $\omega_1^2\omega_3$  are non zero.
Since $H$ belongs to $\cV$, there exists a non-degenerate relation $q=\sum_{ij}q_{ij}\omega_i\omega_j=0$. 
Taking the products with $\omega_k$ for $k=0,\ldots ,3$ we get $q_{00}=q_{02}=q_{11}=q_{13}=0$. 
So the relation $q$ must be a combination of $\omega_0\omega_1, \omega_0\omega_3, \omega_1\omega_2, 
\omega_2^2, \omega_2\omega_3, \omega_3$. 

The pencil $\langle \omega_2, \omega_3\rangle$ is made of forms of constant rank four, 
and is therefore of one of the two possible types found in \cite{mm}:
$$\omega_2 = e_{02}+e_{13}, \qquad \omega_3 = e_{04}+e_{15}$$
for the generic type, will the special type is 
$$\omega_2 = e_{02}+e_{13}, \qquad \omega_3 = e_{03}+e_{14}.$$ 
Note that the pencil $\langle \omega_2, \omega_3\rangle$ spans the singular locus of the 
characteristic surface. In particular it is uniquely defined by $H$ and we call it
{\it the singular pencil}. 

Let us decompose the two other two-forms as 
\begin{equation}\label{0}
\omega_0 = a e_{01}+e_0\wedge f_0+e_1\wedge f_1+\mu, \qquad 
\omega_1 = b e_{01}+e_0\wedge g_0+e_1\wedge g_1+\nu, 
\end{equation}
where $f_0, f_1, g_0, g_1, \mu, \nu$ do not involve $e_0, e_1$. 

\begin{lemma}
The singular pencil $\langle \omega_2, \omega_3\rangle$ must be special.
\end{lemma}

\proof We proceed by contradiction. Suppose the singular pencil is of generic type, and 
choose an adapted basis as above. 

If $q_{01}\ne 0$, we can modify $\omega_0$ and $\omega_1$ by suitable linear combinations of 
$\omega_2$ and $\omega_3$, in such a way that the relation takes the form $q=q_{01}\omega_0\omega_1
+Q(\omega_2, \omega_3)$. (The equation of the characteristic surface itself will not change.) 
In particular, we get the relations
\begin{equation}\label{1}
f_0\wedge\nu+g_0\wedge\mu=f_1\wedge\nu+g_1\wedge\mu=0, \qquad \mu\wedge\nu=0.
\end{equation}

Suppose that the vectors $f_0, f_1, g_0, g_1$ 
are linearly dependant. Then there is a three-plane $P$ 
such that $f_0, f_1, g_0, g_1$ belong to $P$ and $\mu, \nu$ belong to $\wedge^2P$.
But then $\omega_0\wedge\omega_1=e_{01}\wedge\theta$ for some $\theta$ in $\wedge^2P$, 
which has therefore rank at most two. Now the linear relation $q$ yields 
$$q_{01}\theta = 2q_{22}e_{23}+2q_{23}(e_{25}-e_{34})+2q_{33}e_{45}.$$
The right hand side is a form of rank at most two when $q_{23}^2=q_{22}q_{33}$. 
But then $q$ itself has q-rank at most three, a contradiction. 

Suppose now that the vectors $f_0, f_1, g_0, g_1$ are linearly independant. Then the
first two equations in (\ref{1}) are verified if and only if we can write $\mu, \nu$ 
in the form
$$\mu = x f_0\wedge f_1 +y(f_0\wedge g_1+g_0\wedge f_1)+ z g_0\wedge g_1, $$
$$\nu = y f_0\wedge f_1 +z(f_0\wedge g_1+g_0\wedge f_1)+ t g_0\wedge g_1.$$ 
Then the relation $\mu\wedge\nu=0$ amounts to $xt=yz$. Moreover the conditions that 
$\omega_0^3=\omega_1^3=0$ amount to $z=a(xz-y^2)$ and $y=b(yt-z^2)$. If $yz\ne 0$, 
then $ab\ne 0$ and we can write 
$$\omega_0 = \frac{1}{a}(ae_0-f_1)\wedge (ae_1+f_0)+\frac{1}{z}(yf_0+zg_0)\wedge (yf_1+zg_1), $$
$$\omega_1 = \frac{1}{b}(be_0-g_1)\wedge (be_1+g_0)+\frac{1}{y}(yf_0+zg_0)\wedge (yf_1+zg_1). $$
Then the conditions $\omega_0^2\omega_1=0$ and $\omega_0\omega_1^2=0$ reduce to $ay+bz=0$. 
In this case, both triples $\langle ae_0-f_1, be_0-g_1, yf_1+zg_1\rangle$ and 
$\langle ae_1+f_0, be_1+g_0, yf_0+zg_0\langle$ are linearly dependent. This implies that this
six vectors generate a space $Q\subset V_6$ of dimension at most four. But then since 
$\omega_0$ and $\omega_1$ belong to $\wedge^2Q$, the pencil they generate will necessarily 
contain a two-form of rank two, which is a contradiction. 

We also need to consider the case where for example $y=0$, hence also $bz=0$.
If $y=0$ and $b=0$, then the condition 
$\omega_0\omega_1^2=0$ amounts to $x=0$, and then $z=a(xz-y^2)=0$. 
If $y=z=0$, the conditions that $\omega_0^2\omega_1=0$ and $\omega_0\omega_1^2=0$ 
imply that $x=t=0$. But then all the two-forms in $H$ are of the form $e_0\wedge h_0+
e_1\wedge h_1$, hence they have rank at most four and the rank will necessary drop
somewhere, again a contradiction. 

\smallskip
There remains to consider the case where $q_{01}=0$.  In that case we need $q_{03}q_{12}$ to be 
non zero for the $q$ to remain non degenerate, and we can suppose that $q_{03}=q_{12}=1$. 
Then the relation $q$ implies that 
\begin{equation}\label{2}
\mu\wedge e_4+\nu\wedge e_2 = \mu\wedge e_5+\nu\wedge e_3 = 0.
\end{equation}
On the other hand, the fact that the products of $\omega_0$, $\omega_1$ with $\omega_2^2$, 
$\omega_2\omega_3$, $\omega_3^2$ vanish imply that we can write $\mu$ and $\nu$ 
in the form 
$$\mu = xe_{24}+y(e_{25}+e_{34})+ze_{35}, \qquad \nu = x'e_{24}+y'(e_{25}+e_{34})+z'e_{35}.$$
But then $(\ref{2})$ readily imply that $\mu=\nu=0$. We are then again in a situation where 
all the two-forms in $H$ are of the form $e_0\wedge h_0+e_1\wedge h_1$, which is not possible. 
\qed

\medskip We can therefore suppose that the singular pencil $\langle \omega_2, \omega_3\rangle$
is of special type, as given above. Then $\omega_2^2$, $\omega_2\omega_3$, $\omega_3^2$
give $e_{0123}, e_{0124}, e_{0134}$.  Let us decompose $\omega_0$ and $\omega_1$ as in (\ref{0}). 
The fact that their products with $\omega_2^2$, $\omega_2\omega_3$, $\omega_3^2$
vanish means that $\mu$ and $\nu$ do not involve $e_5$. 

\smallskip
Suppose that $q_{01}\ne 0$, so that the relation $q$ can be normalized as before to the form
$q=\omega_0\omega_1+Q(\omega_2,\omega_3)$. We deduce once again that 
$f_0\wedge\nu+g_0\wedge\mu=f_1\wedge\nu+g_1\wedge\mu=0$, and moreover that $f_1\wedge g_0-f_0\wedge g_1$ does not involve $e_5$. In particular we get the relations
\begin{equation}\label{3}
f_{05}\nu+g_{05}\mu=f_{15}\nu+g_{15}\mu=0, \qquad f_{15}g_0-g_{05}f_1-f_{05}g_1+g_{15}f_0=0.
\end{equation}
Note that $\mu$ and $\nu$ cannot both be zero, since otherwise  all the two-forms in $H$ are of the form 
$e_0\wedge h_0+e_1\wedge h_1$, which is not possible. Hence the first two relations imply that 
$f_{05}g_{15}=f_{15}g_{05}$, and the third one, multiplied by $f_{05}$, yields
$$f_{15}(g_{05}f_0+f_{05}g_0)=f_{05}(g_{05}f_1+f_{05}g_1).$$
But then the two-form $g_{05}\omega_0+f_{05}\omega_1$ can be expressed only in terms of $e_0, e_1$ 
and $g_{05}f_0+f_{05}g_0$ (is $f_{05}\ne 0$, or $g_{05}f_1+f_{05}g_1$ if $f_{15}\ne 0$), hence has rank  at most two: a contradiction!

If $f_{05}=f_{15}=0$, then $g_{05}$ and $g_{15}$ cannot both be zero, for otherwise the 
two-forms in $H$ would not involve $e_5$ at all, and thus the rank would be at most four on $H$ and 
would have to drop somewhere. Then the equations (\ref{3}) imply that $\mu=0$ and that  $f_0$ and
$f_1$ are proportional. But then $\omega_0$ can be written in terms of three vectors only, and thus
has rank two, a contradiction again. 

\smallskip
Finally, suppose that $q_{01}=0$, and then normalize as before the relation $q$ to 
$q=\omega_0\omega_3+\omega_1\omega_2+Q(\omega_2,\omega_3)$, which is of maximal q-rank for
any $Q$. The existence of such a relation is equivalent to the equations 
\begin{equation}\label{4}
\mu\wedge e_3+\nu\wedge e_2 = \mu\wedge e_4+\nu\wedge e_3 = 0, \quad f_{05}=g_{15}=f_{15}-g_{05}=0.
\end{equation}
The fact that the pencil $\langle \omega_0, \omega_1\rangle$ contains only forms of rank
four is equivalent to the following conditions:
\begin{equation}\label{5}
\mu\wedge f_0=\nu\wedge g_1=0, \quad  \nu\wedge f_0=\mu\wedge (f_1-g_0), \quad 
\mu\wedge g_1=\nu\wedge (g_0-f_1).
\end{equation}
Moreover we need $\omega_0^2\omega_3=\omega_1^2\omega_2=0$, which amount to 
$\mu\wedge_3=\nu\wedge e_3=0$, and $\omega_0\omega_1\omega_2=\omega_0\omega_1\omega_3=0$, 
which amount to $\mu\wedge e_4=\nu\wedge e_2=0$. Therefore
there exists two scalars $m, n$ such that $\mu = me_{34}$ and $\nu = n e_{23}$,
and the first two equations of (\ref{4}) are then verified. Moreover $\omega_0^2\omega_2=\omega_1^2\omega_3\ne 0$
if and only if $m=-n\ne 0$. The conditions (\ref{5}) then reduce to 
\begin{equation}\label{6}
g_{14}=f_{02}=0, \qquad f_{04}=g_{02}-f_{12}, \qquad g_{12}=f_{14}-g_{04}.
\end{equation}
We also need that $f_{15}=g_{05}\ne 0$ (for otherwise $e_5$ is not involved at all and the
rank has to drop somewhere). If those conditions are satisfied, then the characteristic
surface has the correct equation and the relation $q$ is verified for some $Q$. The last
condition to check, in order to get a four plane $H$ that belongs to $\cV$, is 
that this is the only relation in $Sym^2H$. 

\begin{lemma} 
There is no other relation than $q$. 
\end{lemma}

\proof We grade $\wedge^4V_6$ by the degrees on $e_0, e_1$, respectively 
$e_2, e_3, e_4$, and $e_5$. 
For example $\omega_2^2$, $\omega_2\omega_3$, $\omega_3^2$ span the space of forms of multidegree
$(2,2,0)$. We note that among the binomials in the $\omega_i$, the only ones having terms of
multidegree $(1,2,1)$ are $\omega_0^2$ (which gives $e_{1345}$), $\omega_1^2$ (which gives 
$e_{0235}$), and $\omega_0\omega_1$ (which gives $e_ {0345}+e_{1235}$); those three 
terms are independent. Moreover, among the remaining monomials, only $\omega_0\omega_3$, $\omega_1\omega_2$, $\omega_0\omega_2$, $\omega_1\omega_3$ have terms of multidegree $(2,1,1)$ ($e_{0135}$ for the first two, 
$e_{0125}$ for the third one, $e_{0145}$ for the last one). This readily implies that the dimension of $H\wedge H$ is nine. 
\qed

\medskip
Multiplying $e_0, e_1$ by $m$ and rescaling, we may suppose that $m=1$. 
We get the following result:

\begin{theorem}
Suppose that the characteristic surface $S_H$ of a four plane $H$ in $\cV$ 
is Cayley's ruled cubic surface. Then there exists a basis of $V_6$, and coefficients
$A,B,C,D,E$, such that
$$\omega_0 = e_{15}+e_{34}+Ae_{03}+2Be_{04}-Be_{12}+Ce_{13}+De_{14},$$
$$\omega_1 = e_{05}-e_{23}+Be_{02}+Ee_{03}-De_{04}+2De_{12}+Ee_{13},$$
$$\omega_2 = e_{02}+e_{13}, \qquad \omega_3 = e_{03}+e_{14}.$$
Conversely, any such four plane belongs to $\cV$ and its characteristic 
surface is Cayley's ruled  surface. 
\end{theorem}

In particular, letting $A=B=C=D=E=0$ we get the following four plane of skew-symmetric matrices:
$$\begin{pmatrix}
 0& 0& c& d& 0& b\\
 0& 0& 0& c& d & a \\ 
-c& 0& 0& -b& 0& 0\\
-d& -c& b&  0& a& 0\\
0& -d& 0&  -a& 0& 0\\
-b& -a& 0&  0& 0& 0
\end{pmatrix}$$
The dependance relation is simply $q=\omega_0\omega_3+\omega_1\omega_2=0.$


\subsection{Cones} 
Our final result in this section shows that the characteristic map is not surjective:

\begin{proposition}
Let $H\in\cV$ be such that the characteristic surface $S_H$ is a cone. 
Then $S_H$ is reducible. 
\end{proposition}

\proof Suppose that $S_H$ is a cone over $\omega_0$. In particular $\omega_0^3=0$, so $\omega_0$ 
has rank four and there is a unique $V_4\subset V_6$ such that $\omega_0$ belongs to $\wedge^2V_4$. 
Moreover $\omega_0^2$ is a generator of $\wedge^4V_4$, and the fact that $\omega_0^2\omega=0$ for 
all $\omega\in H$ is equivalent to the fact that $H\subset V_4\wedge V_6$. More concretely, we can 
complete $V_4$ with two independent vectors $e_5, e_6$ such that every $\omega\in H$ can be written
in the form 
$$\omega = e_5\wedge\alpha +e_6\wedge \beta +\theta, $$
where $\alpha, \beta$ belong to $V_4$ and $\theta$ to $\wedge^2V_4$.

\begin{lemma}
The characteristic surface $S_H$ is a cone if and only if $\omega_0$ is $H$-symmetric,
in the sense that $\omega_0\wedge\alpha\wedge\beta=0$ for all $\omega = e_5\wedge\alpha +e_6\wedge \beta +\theta\in H$. 
\end{lemma}

A convenient way to use this lemma is to denote $\tilde\omega_0(\alpha,\beta)=\omega_0\wedge\alpha\wedge\beta$,
where $\tilde\omega_0$ is now a non degenerate skew-symmetric form on $V_4$ (although defined only up to scalar). 

\medskip

\section{The reducible case} 
In this section we discuss the case where the characteristic surface $S_H$ is a reducible
cubic. Equivalently, $S_H$ contains a plane $P$. Since $H$ belongs to $\cV$, this plane must 
be made of two-forms of constant rank four. Such planes have been classified in \cite{mm}; 
up to the action of $PGL_6$ there are exactly four different types, which we consider
case by case.

\medskip\noindent {\it Type 1}. The first type consists in three-planes contained in $\wedge^2\CC^5$; it is represented by the plane $L$ generated by 
$$\begin{array}{rcl}
\omega_1 & = & e_1\wedge e_4+e_2\wedge e_3, \\
\omega_2 & = & e_1\wedge e_5+e_2\wedge e_4, \\
\omega_3 & = & e_2\wedge e_5+e_3\wedge e_4.
\end{array}$$
Note that $Sym^2L\subset\wedge^4\CC^5\subset\wedge^4\CC^6$, and therefore any tensor in 
$L\wedge L$ has rank two. Suppose there exists $\lambda\in L$ and $\omega\in\wedge^2\CC^6$ 
such that $\omega\wedge\lambda$ belongs to $L\wedge L$. Then necessarily 
$\omega$ belongs to $\wedge^2\CC^5$. But then $H$ must be a 
four plane in$\wedge^2\CC^5$, in which the variety of tensors 
of rank at most two has codimension $3$. So $H$ necessarily 
contains some tensors of rank two, a contradiction. 

Now suppose that $\omega\wedge\omega$ belongs to $L\wedge L$. 
Let us decompose $\omega=\phi+\alpha\wedge e_6$. Then 
$\omega\wedge\omega=\phi\wedge\phi+2\phi\wedge\alpha\wedge e_6.$
In particular we need $\phi\wedge\alpha=0$. If $\alpha=0$, then
$\omega$ belongs to $\wedge^2\CC^5$ and we get a contradiction
as before. Otherwise, $\alpha$ must divide $\phi$ and then 
$\omega$ has rank two, a contradiction again. 
We conclude that no three plane of type 1 can be contained in
a four plane of $\cV$. 

\medskip\noindent {\it Type 2}. The second type is provided by the plane $L$ generated by 
$$\begin{array}{rcl}
\omega_1 & = & e_0\wedge e_4-e_1\wedge e_3, \\
\omega_2 & = & e_0\wedge e_5-e_2\wedge e_3, \\
\omega_3 & = & e_1\wedge e_5-e_2\wedge e_4, 
\end{array}$$
An easy computation shows that the $L$ is stabilized by a subgroup
$Stab_L$ of $GL_6$ isomorphic to $\CC^*\times GL_3$. Moreover the action of $Stab_L$ on $L$ is equivalent to the action of  
$\CC^*\times GL_3$ on $A_3(\CC)$, the space of $3\times 3$ skew-symmetric matrices, given by $(z,A): X\mapsto zA^tXA$. 
In particular the action of $Stab_L$ on $L$ has exactly two orbits, the origin and its complement; and the action of 
$Stab_L$ on $Sym^2L$ has exactly four orbits, given by the q-rank. 

For $H$ to belong to $\cV$, we need $\omega\wedge\omega$ to be an element of 
$L\wedge L$ of q-rank three. We may suppose this element, up 
to the action of $Stab_L$, to be 
$$-\frac{1}{2}(\omega_1^2+\omega_2^2+\omega_3^2)=e_{0134}+e_{0235}+e_{1245}.$$
Then we recover $\omega$ by applying $pf^*$, which yields (up to scalar) 
$$\omega = e_0\wedge e_3+e_1\wedge e_4+e_2\wedge e_5.$$
We also need to consider the case where there exists $\lambda\in L$ such that $\omega\wedge\lambda$ belongs to $L\wedge L$. 
Up to the action of $Stab_L$ we may suppose that $\lambda=\omega_1$. 
Then a straightforward computation shows that (modulo $L$) 
$\omega$ must be a combination of $e_0\wedge e_1$, $e_0\wedge e_3$, $e_1\wedge e_4$, $e_3\wedge e_4$.  In particular $\omega$ 
belongs to $\wedge^2M$, where $M=\langle e_0,e_1,e_3,e_4\rangle$. Therefore the line joining $\omega$ to $\omega_1$, which
also belongs to $\wedge^2M$, will necessary contain an element of rank two, and we get a 
contradiction. 

The conclusion of this discussion is that up to equivalence,
a three plane $L$ of type two can be uniquely extended to a 
four plane $H$ in $\cV$, which in a suitable basis is the four 
plane of matrices of the form:
$$\begin{pmatrix}
 0& 0& 0& d& a& b\\
 0& 0& 0& -a& d& c \\ 
0& 0& 0& -b& -c& d\\
d& a& b&  0& 0& 0\\
-a& d& c&  0& 0& 0\\
-b& -c& d&  0& 0& 0
\end{pmatrix}$$
The associated cubic surface has equation $d(a^2+b^2+c^2+d^2)=0$. It is the union of the plane $L$ and of a smooth quadric.  

\medskip\noindent {\it Type 3}. The third type is provided by the plane $L$ generated by 
$$\begin{array}{rcl}
\omega_1 & = & e_0\wedge e_2+e_1\wedge e_3, \\
\omega_2 & = & e_0\wedge e_3+e_1\wedge e_4, \\
\omega_3 & = & e_0\wedge e_4+e_1\wedge e_5, 
\end{array}$$
In this case, $L\wedge L\simeq Sym^2L$ is the space of four-forms divisible by
$e_0\wedge e_1$. In particular $L\wedge L$ does not contain any rank six element, 
so our general strategy does not apply. 

Instead, let us try directly to complete $L$ into a four-plane $H\in\cV$, with some two-form $\omega$ such 
that $\omega\wedge\omega$ belongs to $L\wedge L$. 
We can decompose $\omega=e_0\wedge u_0
+e_1\wedge u_1+\psi$, where $\psi$ does not involve $e_0, e_1$. Then $\omega\wedge\omega$
is divisible by $e_0\wedge e_1$ if and only if 
$$u_0\wedge\psi = u_1\wedge\psi = 0 \quad \mathrm{and} \quad \psi\wedge\psi = 0.$$
This means that $\psi$ has rank two ($\psi$ cannot be zero, for otherwise $H\wedge H=
L\wedge L$) and is divisible by $u_0$ and $u_1$. Note that $u_0, u_1$ cannot be dependent,
since otherwise $\omega$ would have rank two. So $\psi$ must be a multiple of $u_0\wedge u_1$ 
and after a suitable normalization we can suppose that 
$$\omega = e_0\wedge u_0+e_1\wedge u_1+u_0\wedge u_1.$$
Note that an element in $L$ is of the form $\lambda = e_0\wedge \ell_0+e_1\wedge \ell_1$,
and the planes $\langle \ell_0,\ell_1\rangle$ span a Veronese surface in $G(2,5)$. It is 
easy to see that $(\omega\wedge L)\cap (L\wedge L)=0$ if and only if the plane 
$\langle u_0,u_1\rangle$ does not belong to this surface. This ensures that the kernel
of the map $Sym^2H\ra H\wedge H$ is one-dimensional. For $H$ to belong to $\cV$, we 
finally need that a generator of this kernel be of maximal q-rank, which is easily checked to
be true in general: take for example $u_0=e_5$ and $u_1=e_2$, then the relation is 
$$\omega^2+\omega_1\omega_3+\omega_2^2=0,$$ 
which is non degenerate. The four-plane of skew-symmetric matrices is:
$$\begin{pmatrix}
 0& 0& b& c& d& a\\
 0& 0& a& b& c& d \\ 
-b& -a& 0& 0& 0& -a\\
-c& -b& 0&  0& 0& 0\\
-d& -c& 0&  0& 0& 0\\
-a& -d& a&  0& 0& 0
\end{pmatrix}$$

A general point in $H$ is of the form $\phi=e_0\wedge (su_0+\ell_0)+e_1\wedge (su_1+\ell_1)+
su_0\wedge u_1$, and has rank four or less when $s \ell_0\wedge \ell_1\wedge u_0\wedge u_1=0$. 
This implies that the characteristic surface is the union of  $L$ and a cone over a conic 
in $L$. This conic may be singular: consider for example the case where $u_0=e_5$ and 
$u_1=e_3$; we get the following four plane of skew-symmetric matrices:
$$\begin{pmatrix}
 0& 0& b& c& d& a\\
 0& 0& 0& a+b& c& d \\ 
-b& 0& 0& 0& 0& 0\\
-c& -a-b& 0&  0& 0& -a\\
-d& -c& 0&  0& 0& 0\\
-a& -d& 0&  a& 0& 0
\end{pmatrix}$$

\medskip\noindent {\it Type 4}. The fourth type is represented by the plane $L$ generated by 
$$\begin{array}{rcl}
\omega_1 & = & e_0\wedge e_3+e_1\wedge e_2, \\
\omega_2 & = & e_0\wedge e_4+e_2\wedge e_3, \\
\omega_3 & = & e_0\wedge e_5+e_1\wedge e_3.
\end{array}$$
Then $L\wedge L\simeq Sym^2L$ is generated by $f_{14}+f_{25}, f_{24}, f_{34}, f_{15}, f_{35}, f_{45}$. In particular
an element of $L\wedge L$ can be written as $f_4\wedge g_4+f_5\wedge g_5$, hence never has rank six. 

So let us directly try to complete $L$ into a four-plane $H\in\cV$ with a two-form $\omega$ such 
that $\omega\wedge\omega$ belongs to $L\wedge L$. Let us decompose 
$$\omega = \psi + u_4\wedge e_4+u_5\wedge e_5+ze_4\wedge e_5,$$
where $\psi, u_4, u_5$ do not involve $e_4, e_5$. Then we have 
$$\frac{1}{2}\omega\wedge\omega =\frac{1}{2}\psi\wedge\psi +\psi\wedge u_4\wedge e_4+\psi\wedge u_5\wedge e_5
+(z\psi-u_4\wedge u_5)\wedge e_4\wedge e_5.$$
For this to belong to $L\wedge L$, we need that $z\psi=u_4\wedge u_5$. If $z\ne 0$, we simply get 
that $\omega\wedge\omega =0$, which is excluded. If $z=0$, then $u_4$ and $u_5$ must be colinear, 
and after a change of basis we may suppose that $u_5=0$. Then we need $\psi\wedge u_4$ to be a 
combination of $e_{012}$ and $e_{023}$, in which case $\psi\wedge u_4\wedge e_4$ is a linear combination
of $\omega_1\wedge\omega_2$ and $\omega_2\wedge\omega_2$. But then the kernel of $Sym^2H\lra H\wedge H$ is 
generated by a quadratic relation of the form $\omega^2=Q(\omega_1,\omega_2)$, which contradicts 
the non degeneracy condition. 

There remains the possibility that $(\omega\wedge L)\cap (L\wedge L)\ne 0$. Another computation leads to the same conclusion. 

\medskip 
We summarize this discussion as follows:

\begin{theorem} Let $H\in\cV$ be such that the characteristic surface $S_H$ contains 
a plane $L$. Then either:
\begin{enumerate}
\item $L$ is of type 2 and the residual component of $S_H$ is a smooth quadric;
\item $L$ is of type 3 and the residual component of $S_H$ is a quadratic cone, whose vertex 
is outside $P$.
\end{enumerate}
\end{theorem}

\bigskip

\noindent {\sc E.V. Ferapontov}, 
Department of Mathematical Sciences, Loughborough University,
Loughborough, Leicestershire, LE11 3TU, UK. 

\noindent\texttt{e.v.ferapontov@lboro.ac.uk}

\smallskip\noindent {\sc L. Manivel}, Institut de Math\'ematiques de Toulouse, UMR 5219, Universit\'e de Toulouse, CNRS, UPS IMT F-31062 Toulouse Cedex 9, France.

\noindent\texttt{manivel{@}math.cnrs.fr}
\end{document}